\documentclass[graybox]{svmult}
\usepackage{amsmath,amssymb}
\usepackage{mathptmx}       
\usepackage{helvet}         
\usepackage{courier}        
\usepackage{type1cm}        
\usepackage{makeidx}         
\usepackage{graphicx}        
\usepackage{multicol}        
\usepackage[bottom]{footmisc}%
\makeindex             
\begin{document}
\author{Susanne C. Brenner, Christopher B. Davis, and Li-yeng Sung}
\title*{Additive Schwarz preconditioners for a state constrained elliptic
distributed  optimal control problem
 discretized by a partition of unity method}
\titlerunning{AS preconditioners for a state constrained elliptic optimal control problem}
\maketitle
\abstract*{
We present additive Schwarz preconditioners for a class of elliptic optimal control problems
discretized by a partition of unity method.  The discrete problem is solved by
   a primal-dual active set algorithm, where the auxiliary system in each iteration is solved
 by a preconditioned conjugate gradient method based on additive Schwarz preconditioners.
Condition number estimates are given and verified by a numerical example.
}
\section{Introduction}\label{sec:Introduction}
In this work, we are interested in solving a model elliptic optimal control problem of
the following form:  Find $(y,u)\in H^1_0(\Omega)\times L_2(\Omega)$
 that minimize the functional
\begin{displaymath}
J(y,u)=\frac{1}{2}\int_\Omega (y-f)^2 dx + \frac{\beta}{2}\int_\Omega u^2 dx
\end{displaymath}
subject to
\begin{equation}
-\Delta y=u \ \text{in} \ \Omega,  \quad
y = 0  \ \text{in} \ \partial \Omega,
\label{eq:Poisson}
\end{equation}
and
 $y \leq \psi \ \text{in} \ \Omega$,
where $\Omega$ is a convex polygon in $\mathbb{R}^2$ and $f\in L_2(\Omega).$  We also assume
$\psi\in C^2(\Omega)\cap H^3(\Omega)$ and $ \psi > 0   \ \text{on}\ \partial\Omega.$
\par
 Due to the elliptic regularity (cf. \cite{Grisvard:1985:EPN}) for \eqref{eq:Poisson},
 we can reformulate the model problem as follows:
Find $y\in K$ such that
\begin{equation}\label{eq:cts_prob}
	y=\mathop{\rm argmin}_{v\in K}	\left[\frac{1}{2}a(v,v)-(f,v)\right],
\end{equation}
 where
 $K = \{v\in H^2(\Omega)\cap H^1_0(\Omega):v \leq \psi \ \text{in} \ \Omega\}$,
\begin{equation*}
		a(w,v) =\beta\int_{\Omega} \Delta w \Delta v dx + \int_{\Omega} wv dx \quad\text{and}\quad
		(f,v) =\int_{\Omega} fv dx.
\end{equation*}  Once $y$ is calculated, then $u$ can be determined by $u=-\Delta y.$
\par
 The minimization problem \eqref{eq:cts_prob}
  was discretized in \cite{BDS_PUM_VI} by a partition of unity method (PUM).
 The goal of this paper is to use the ideas in \cite{BDS_ASM_VI} for an obstacle problem
  of clamped Kirchhoff plates to develop preconditioners for
 the discrete problems  in \cite{BDS_PUM_VI}.
 We refer to these references for technical details and only present the important results here.
%
\section{The Discrete Problem}\label{sec:Discrete}
 We will use a variant of the PUM  (cf. \cite{MB_PUFEM,GS_PPUM,BBO_MESHLESS,OKH_PUM})
  to construct a conforming
 approximation space $V_h\subset H^2(\Omega)\cap H^1_0(\Omega)$.
 Below we present an overview of the construction of $V_h$.
\par
 We take $V_h$ to be
 $\sum_{i=1}^n \phi_i V_i$, where $\{\phi_i\}_{i=1}^n$ is a set of $C^1$
 piecewise polynomial flat-top partition of unity functions over $\Omega$ and $\{V_i\}_{i=1}^n$ are the local approximation spaces.  We denote the patch $\Omega_i$ as the support of $\phi_i$
  and the flat-top part of $\Omega_i$ is $\Omega_i^{\text{flat}}=\{x\in \Omega : \phi_i(x) = 1\}$
   (for an example see Figure 2.1 in \cite{BDS_ASM_VI}).  Each
   $V_i$
    consists of biquadratic polynomials satisfying the Dirichlet boundary condition
     of \eqref{eq:Poisson}.  Basis functions for $V_i$ are tensor product Lagrange
     polynomials and the interpolation nodes are
 distributed uniformly over $\Omega_i^{\text{flat}}.$  This allows us to select
 basis functions for $V_h$ that satisfy the Kronecker delta property.
\par
Let $\mathcal{N}_h$ be the set of all interior interpolation nodes used in
 the construction of $V_h.$ The discrete problem is to find $y_h\in K_h$ such that
\begin{equation}
y_h=\mathop{\rm argmin}_{v\in K_h}
\left[\frac{1}{2}a(v,v)-(f,v)\right],
\label{eq:discrete_prob}
\end{equation}
where
  $K_h = \{v\in V_h: v(p) \leq \psi(p) \ \forall p \in \mathcal{N}_h\}$.
\par
 By introducing a Lagrange multiplier $\lambda_h:\mathcal{N}_h \to \mathbb{R}$, the minimization problem
  \eqref{eq:discrete_prob} can be rewritten in  the following form: Find $y_h \in K_h$ such that
\begin{alignat*}{3}
a(y_h,v)-(f,v) &= -\sum_{p\in \mathcal{N}_h} \lambda_h(p)v(p) &\qquad& \forall\, v\in V_h, \\
\lambda_h(p) &= \max(0,\lambda_h(p)+c(y_h(p)-\psi(p))) &\qquad& \forall\, p \in \mathcal{N}_h,
\end{alignat*}
 where $c$ is a (large) positive number.
This system can then be solved by a primal-dual active set (PDAS) algorithm
(cf. \cite{BIK_PD,BK_PD,HIK_PDAS,IK_LM}).  Given the $k$-th approximation $(y_k, \lambda_k)$,
the $(k+1)$-st iteration of the PDAS algorithm is to find $(y_{k+1},\lambda_{k+1})$ such that
\begin{alignat}{3}\label{eq:linearsystem}
a(y_{k+1},v)-(f,v) &= -\sum_{p\in \mathcal{N}_h} \lambda_{k+1}(p)v(p) &\qquad& \forall v\in V_h,
\notag \\
y_{k+1}(p) &= \psi(p) &\qquad& \forall p \in \mathfrak{A}_k,\\
\lambda_{k+1}(p) &= 0 &\qquad& \forall p \in \mathcal{N}_h \backslash \mathfrak{A}_k,\notag
\end{alignat}
where $\mathfrak{A}_k = \{ p\in \mathcal{N}_h : \lambda_k(p) + c(y_k(p)-\psi(p)) > 0  \}$
 is the set of active nodes determined from the approximations $(y_k,\lambda_k).$
 Below we present preconditioners for the linear systems encountered in \eqref{eq:linearsystem}.
\section{The Preconditioners}
\label{sec:pu}
The additive Schwarz preconditioners (cf. \cite{DW:1987:AS})
 will be applied to a system associated with a subset $\tilde{\mathcal{N}}_h$
 of $\mathcal{N}_h$.
  Let   $\tilde{T}_h: V_h \to V_h$ be defined by

\begin{equation*}
(\tilde{T}_h v)(p) = \left\{\begin{array}{cc}v(p) & \text{if} \ p\in
\tilde{\mathcal{N}}_h \\ 0 & \text{if} \ p \not\in \tilde{\mathcal{N}}_h \end{array}\right..
\end{equation*}
The approximation space for the subproblem is  $\tilde{V}_h = \tilde{T}_h V_h$.
 The associated stiffness matrix is a symmetric positive definite operator
 $\tilde{A}_h:\tilde{V}_h \to \tilde{V}_h'$ defined by
\begin{equation*}
\langle \tilde{A}_h v, w \rangle = a(v,w) \qquad \forall\, v,w\in \tilde{V}_h,
\end{equation*}
where $\langle \cdot , \cdot \rangle $ is the canonical bilinear form on
 $\tilde{V}_h' \times \tilde{V}_h.$
\par\bigskip\noindent
{\bf A One-Level Method}\quad
Here we introduce a collection of shape regular
subdomains $\{D_j\}_{j=1}^J$ with $\text{diam} \ D_j \approx H$ that overlap
 with each other by at most $\delta.$  Associated with each subdomain is a
  function space $V_j\subset \tilde{V}_h$ whose members vanish at the nodes outside $D_j.$
Let $A_j:V_j\to V_j'$ be defined by
\begin{displaymath}
\langle A_j v, w \rangle = a(v,w) \qquad \forall v,w\in V_j.
\end{displaymath}
The one-level additive Schwarz preconditioner $B_{\text{OL}}:V_h'\to V_h$ is defined by
\begin{equation*}
B_{\text{OL}} = \sum_{j=1}^J I_j A_j^{-1} I_j^t,
\end{equation*}
where $I_j : V_j \to \tilde{V}_h$ is the natural injection.

Following the arguments  in \cite{BDS_ASM_VI},  we can obtain the following theorem.
\begin{theorem}\label{thm:One_level}
	There exists a positive constant $C_{\text{OL}}$ independent of $H$, $h$, $J$, $\delta$ and $\tilde{N}_h$
   such that
	\begin{displaymath}
	\kappa(B_{\text{OL}}\tilde{A}_h) \leq C_{\text{OL}}\delta^{-3}H^{-1}.
	\end{displaymath}
	\end{theorem}
\par\bigskip\noindent
{\bf A Two-Level Method} \quad
Let $V_H\subset H^2(\Omega)\cap H^1_0(\Omega)$ be a coarse
approximation space based on the construction in Section \ref{sec:Discrete}
where $H>h.$ We assume the patches of $V_H$ are of comparable size to the subdomains
$\{D_j\}_{j=1}^J.$  Let $\Pi_h:H^2(\Omega)\cap H^1_0(\Omega)\to V_h$ be the
nodal interpolation operator.
   We define $V_0\subset \tilde{V}_h$ by
$V_0=T_h\Pi_h V_H$,
 and $A_0:V_0\to V_0'$ by
\begin{equation*}
\langle A_0 v,w\rangle = a(v,w) \qquad \forall\, v,w\in V_0.
\end{equation*}
The two-level additive Schwarz preconditioner $B_{\text{TL}} : V_h' \to V_h$ is given by
\begin{equation*}
B_\text{TL} = \sum_{j=0}^J I_j A_j^{-1} I_j^t,
\end{equation*}
where $I_0:V_0 \to \tilde{V}_h$ is the natural injection.
Following the arguments in \cite{BDS_ASM_VI}, we can obtain the following theorem.
\begin{theorem}\label{thm:Two_level}
	There exists a positive constant $C_{\text{TL}}$ independent of
  $H$, $h$, $J$, $\delta$ and $\tilde{N}_h$  such that
	\begin{displaymath}
\kappa(B_{\text{TL}}A_h) \leq C_{\text{TL}}
   \min\big((H/h)^4,\delta^{-3}H^{-1}).
	\end{displaymath}
	\end{theorem}
\begin{remark} \label{rem:TL}
 The two-level method is scalable as long as $H/h$ remains bounded.
\end{remark}
\section{A Numerical Example}\label{sec:Numerics}
We consider Example 4.2 in \cite{BDS_PUM_VI}, where
 $\Omega = (-0.5, 0.5)^2$, $\beta=0.1,$ $\psi=0.01$, and
 $f =10(\sin(2\pi(x_1+0.5)) + (x_2+0.5))$.  We discretize
 \eqref{eq:discrete_prob} by the PUM with uniform rectangular patches
so that $h \approx 2^{-\ell}$, where $\ell$ is the refinement level. As $\ell$ increases from 1 to
8, the number of degrees of freedom increases from 16 to 586756. The discrete variational
inequalities are solved by the PDAS
algorithm presented in Section \ref{sec:Discrete}, with $c=10^8$.
\par
For the purpose of comparison, we first solve the auxiliary systems in each iteration of
the PDAS algorithm by the conjugate gradient (CG) method without a preconditioner. The
average condition number during the PDAS iteration and the time to solve the variational
inequality are presented in Table \ref{tab:nopc}. The PDAS iterations fail to stop (DNC) within 48
hours beyond level 6.
\begin{table}[htbp]
\caption{Average condition number ($\kappa$) and time to solve 
 ($t_\text{solve}$) in seconds by the CG algorithm}
\label{tab:nopc}
\begin{tabular}{|c|c|c|}
\hline
$\ell$ & $\kappa$ & $t_\text{solve}$ \\ \hline
1  & 3.1305$\times10^{+2}$ & 2.6111$\times10^{-2}$\\ \hline
2  & 9.1118$\times10^{+3}$ & 1.0793$\times10^{-1}$\\ \hline
3  & 2.0215$\times10^{+5}$ & 9.7842$\times10^{-1}$\\ \hline
4  & 3.3705$\times10^{+6}$ & 3.3911$\times10^{+1}$\\ \hline
5  & 6.4346$\times10^{+7}$ & 6.2173$\times10^{+2}$\\ \hline
6  & 1.0537$\times10^{+9}$ & 8.8975$\times10^{+3}$\\ \hline
7 & DNC & DNC \\ \hline
8 & DNC & DNC \\ \hline
\end{tabular}
\end{table}
\par
We then solve the auxiliary systems by the preconditioned conjugate gradient (PCG)
method, using the additive Schwarz preconditioners associated with $J$ subdomains. The
mesh size $H$ for the coarse  space $V_H$ is $\approx 1/\sqrt{J}$.
 We say the PCG
method has converged if $\|Br\|_2 \leq 10^{-15}\|b\|_2$,
where $B$ is the preconditioner, $r$ is the residual,
and $b$ is the load vector. The initial guess for the PDAS algorithm is taken to be the solution
at the previous level, or 0 if $2^{2\ell} = J$.
To obtain a good initial guess for the two-level method, the one-level method
 is used when $2^{2\ell}=J$. The subdomain problems and the coarse problem are
solved by a direct method based on the Cholesky factorization on independent processors.
\par\bigskip\noindent
{\bf Small Overlap} \quad
Here we apply the preconditioners in such a way that $\delta \approx h$.
 The averaged condition numbers of the linear systems over the PDAS
 iterations are presented in Tables \ref{tab:OLSO_k} and \ref{tab:TLSO_k}.
  We can see that
these condition numbers are significantly smaller than those for the
 unpreconditioned case and
the condition numbers for the two-level method are smaller than those for the one-level method.
For each $\ell$, as $J$ increases the condition numbers for the two-level method are decreasing,
 which demonstrates the scalability of the two-level method (cf. Remark~\ref{rem:TL}).
\begin{table}[htbp]
\caption{One-level small overlap: average condition number}
\begin{tabular}{|c|c|c|c|c|}
\hline
$\ell$ & $J=4$ & $J=16$ & $J=64$ & $J=256$ \\ \hline
 1 & 1.0000$\times10^{+0}$ & - & - & -  \\ \hline
 2 & 4.9489$\times10^{+0}$ & 7.4007$\times10^{+0}$ & - & -  \\ \hline
 3 & 1.5165$\times10^{+1}$ & 4.4136$\times10^{+1}$ & 6.6141$\times10^{+1}$ & -  \\ \hline
 4 & 7.8249$\times10^{+1}$ & 1.9021$\times10^{+2}$ & 5.3590$\times10^{+2}$ & 8.1948$\times10^{+2}$  \\ \hline
 5 & 6.4747$\times10^{+2}$ & 1.6428$\times10^{+3}$ & 3.1727$\times10^{+3}$ & 9.5009$\times10^{+3}$  \\ \hline
 6 & 5.0797$\times10^{+3}$ & 1.3163$\times10^{+4}$ & 2.5880$\times10^{+4}$ & 5.0481$\times10^{+4}$  \\ \hline
 7 & 4.0710$\times10^{+4}$ & 1.0630$\times10^{+5}$ & 2.1060$\times10^{+5}$ & 4.1597$\times10^{+5}$  \\ \hline
 8 & 3.2674$\times10^{+5}$ & 8.5575$\times10^{+5}$ & 1.7014$\times10^{+6}$ & 3.3814$\times10^{+6}$  \\ \hline

\end{tabular}
\label{tab:OLSO_k}
\end{table}
\begin{table}[htbp]
\caption{Two-level small overlap: average condition number}
\begin{tabular}{|c|c|c|c|c|}
\hline
$\ell$ & $J=4$ & $J=16$ & $J=64$ & $J=256$ \\ \hline
 1 & 1.0000$\times10^{+0}$ & - & - & - \\ \hline
 2 & 5.4624$\times10^{+0}$ & 7.4007$\times10^{+0}$ & - & - \\ \hline
 3 & 1.2293$\times10^{+1}$ & 1.1437$\times10^{+1}$ & 6.6141$\times10^{+1}$ & - \\ \hline
 4 & 2.8578$\times10^{+1}$ & 2.7932$\times10^{+1}$ & 1.2645$\times10^{+1}$ &8.1948$\times10^{+2}$ \\ \hline
 5 & 6.2993$\times10^{+1}$ & 9.1991$\times10^{+1}$ & 4.6130$\times10^{+1}$ & 1.9891$\times10^{+1}$ \\ \hline
 6 & 3.6714$\times10^{+2}$ & 3.4816$\times10^{+2}$ & 1.3100$\times10^{+2}$ & 5.7716$\times10^{+1}$ \\ \hline
 7 & 2.7431$\times10^{+3}$ & 2.1182$\times10^{+3}$ & 1.0314$\times10^{+3}$ & 2.8685$\times10^{+2}$ \\ \hline
 8 & 2.1656$\times10^{+4}$ & 1.4827$\times10^{+4}$ & 9.1992$\times10^{+3}$ & 1.8754$\times10^{+3}$ \\ \hline
\end{tabular}
\label{tab:TLSO_k}
\end{table}
\par
The times to solve the problem for each method are presented in
 Tables~\ref{tab:OLSO_t} and \ref{tab:TLSO_t}.  By comparing them with the
 results in Table~\ref{tab:nopc}, we can see that both methods are superior.
   For comparison purposes, the faster time between the two methods is highlighted
   in red for each
    $\ell$ and $J.$  As $h$ decreases and $J$ increases, the two-level method performs
     better than the one-level method.  This agrees with what one would expect
      from Theorems ~\ref{thm:One_level} and \ref{thm:Two_level}.
\begin{table}[htbp]
\caption{One-level small overlap: time to solve in seconds.  Times highlighted in red are
 faster than the ones for the corresponding two-level method.}
\begin{tabular}{|c|c|c|c|c|}
\hline
$\ell$ & $J=4$ & $J=16$ & $J=64$ & $J=256$ \\ \hline
 1 & 1.7824$\times10^{+0}$ & - & - & - \\ \hline
 2 & \color{red}{3.0485$\times10^{-1}$} & 1.5566$\times10^{+1}$ & - & - \\ \hline
 3 & \color{red}{3.8408$\times10^{-1}$} & \color{red}{1.0783$\times10^{+1}$} & 6.0871$\times10^{+1}$ & - \\ \hline
 4 & \color{red}{2.6069$\times10^{+0}$} & \color{red}{4.1818$\times10^{+1}$} & 9.1811$\times10^{+1}$ & 3.5518$\times10^{+2}$ \\ \hline
 5 & \color{red}{2.5704$\times10^{+1}$} & \color{red}{1.1104$\times10^{+2}$} & 1.5399$\times10^{+2}$ & 3.5482$\times10^{+2}$ \\ \hline
 6 & 2.8261$\times10^{+2}$ & 2.6935$\times10^{+2}$ & 4.0033$\times10^{+2}$ & 4.6376$\times10^{+2}$ \\ \hline
 7 & 5.2566$\times10^{+3}$ & 1.9115$\times10^{+3}$ & 1.4825$\times10^{+3}$ & 1.5815$\times10^{+3}$ \\ \hline
 8 & 1.0946$\times10^{+5}$ & 2.9034$\times10^{+4}$ & 1.1631$\times10^{+4}$ & 6.8551$\times10^{+3}$ \\ \hline
\end{tabular}
\label{tab:OLSO_t}
\end{table}
\begin{table}[htbp]
\caption{Two-level small overlap: time to solve in seconds.  Times highlighted in red are
 faster than the ones for the corresponding one-level method.}
\begin{tabular}{|c|c|c|c|c|}
\hline
$\ell$ & $J=4$ & $J=16$ & $J=64$ & $J=256$ \\ \hline
 1 & 1.7824$\times10^{+0}$ & - & - & - \\ \hline
 2 & 1.0694$\times10^{+0}$ & 1.5566$\times10^{+1}$ & - & - \\ \hline
 3 & 1.0889$\times10^{+0}$ & 1.4261$\times10^{+1}$ & 6.0871$\times10^{+1}$ & - \\ \hline
 4 & 5.5186$\times10^{+0}$ & 5.8349$\times10^{+1}$ & \color{red}{7.0962$\times10^{+1}$} & 3.5518$\times10^{+2}$ \\ \hline
 5 & 3.0950$\times10^{+1}$ & 1.1499$\times10^{+2}$ & \color{red}{1.4285$\times10^{+2}$} & \color{red}{1.4650$\times10^{+2}$} \\ \hline
 6 & \color{red}{2.8179$\times10^{+2}$} & \color{red}{2.0602$\times10^{+2}$} & \color{red}{1.6374$\times10^{+2}$} & \color{red}{1.5015$\times10^{+2}$} \\ \hline
 7 & \color{red}{4.4391$\times10^{+3}$} & \color{red}{1.1894$\times10^{+3}$} & \color{red}{4.6832$\times10^{+2}$} & \color{red}{2.9826$\times10^{+2}$} \\ \hline
 8 & \color{red}{9.0540$\times10^{+4}$} & \color{red}{2.0476$\times10^{+4}$} & \color{red}{3.1224$\times10^{+3}$} & \color{red}{8.8092$\times10^{+2}$} \\ \hline
\end{tabular}
\label{tab:TLSO_t}
\end{table}
\par\bigskip\noindent{\bf Generous Overlap}\quad
Here we apply the preconditioners in such a way that $\delta \approx H$.
 When $J=4$ and $J=16$ both methods fail to converge at $\ell=8$ within 48
 hours due to the large size of the local problems.  The averaged condition numbers
  of the linear systems over the PDAS iterations are presented in
   Tables~\ref{tab:OLLO_k} and \ref{tab:TLLO_k}.  They agree with Theorems ~\ref{thm:One_level} and \ref{thm:Two_level}.  We can also see that these condition numbers
are smaller than those in the case of small overlap.
\begin{table}[htbp]
\caption{One-level generous overlap: average condition number}
\begin{tabular}{|c|c|c|c|c|}
\hline
$\ell$ & $J=4$ & $J=16$ & $J=64$ & $J=256$ \\ \hline
 1 & 1.0000$\times10^{+0}$ & - & - & -  \\ \hline
 2 & 1.0000$\times10^{+0}$ & 7.4007$\times10^{+0}$ & - & -  \\ \hline
 3 & 1.0000$\times10^{+0}$ & 7.8491$\times10^{+0}$ & 6.6141$\times10^{+1}$ & -  \\ \hline
 4 & 1.0000$\times10^{+0}$ & 7.5665$\times10^{+0}$ & 8.4735$\times10^{+1}$ & 8.1948$\times10^{+2}$  \\ \hline
 5 & 1.0000$\times10^{+0}$ & 8.2910$\times10^{+0}$ & 9.6722$\times10^{+1}$ & 1.4803$\times10^{+3}$  \\ \hline
 6 & 1.0000$\times10^{+0}$ & 8.3675$\times10^{+0}$ & 9.8624$\times10^{+1}$ & 1.4780$\times10^{+3}$  \\ \hline
 7 & 1.0000$\times10^{+0}$ & 8.4332$\times10^{+0}$ & 1.0019$\times10^{+2}$ & 1.4998$\times10^{+3}$  \\ \hline
 8 & DNC & DNC & 1.0108$\times10^{+2}$ & 1.5161$\times10^{+3}$  \\ \hline
\end{tabular}
\label{tab:OLLO_k}
\end{table}
\begin{table}[htbp]
\caption{Two-level generous overlap: average condition number}
\begin{tabular}{|c|c|c|c|c|}
\hline
$\ell$ & $J=4$ & $J=16$ & $J=64$ & $J=256$ \\ \hline
 1 & 1.0000$\times10^{+0}$ & - & - & - \\ \hline
 2 & 1.2500$\times10^{+0}$ & 7.4007$\times10^{+0}$ & - & - \\ \hline
 3 & 1.2500$\times10^{+0}$ & 6.2713$\times10^{+0}$ & 6.6141$\times10^{+1}$ & - \\ \hline
 4 & 1.2500$\times10^{+0}$ & 6.4760$\times10^{+0}$ & 1.3273$\times10^{+1}$ & 8.1948$\times10^{+2}$ \\ \hline
 5 & 1.2500$\times10^{+0}$ & 7.1544$\times10^{+0}$ & 1.7516$\times10^{+1}$ & 1.7316$\times10^{+1}$ \\ \hline
 6 & 1.2500$\times10^{+0}$ & 7.4536$\times10^{+0}$ & 2.0683$\times10^{+1}$ & 2.0360$\times10^{+1}$ \\ \hline
 7 & 1.2500$\times10^{+0}$ & 7.6360$\times10^{+0}$ & 2.2223$\times10^{+1}$ & 2.5925$\times10^{+1}$ \\ \hline
 8 & DNC & DNC & 2.4425$\times10^{+1}$ & 2.8208$\times10^{+1}$ \\ \hline
\end{tabular}
\label{tab:TLLO_k}
\end{table}
\par
The times to solve the problem for each method are presented in
Tables~\ref{tab:OLLO_t} and \ref{tab:TLLO_t}.
Again both methods are superior to the unpreconditioned method and the scalability of the
two-level method is observed.
\par
We now compare the generous overlap methods with the small overlap methods.
 In Tables~\ref{tab:OLLO_t} and \ref{tab:TLLO_t}, the times in red  are
 the ones where
  the method with generous overlap outperforms the method with small overlap.
  It is evident from Table~\ref{tab:TLLO_t} that the performance of the
   two-level method with generous overlap suffers from a high communication cost for
   small $h$ and large $J$.
\begin{table}[htbp]
\caption{One-level generous overlap: time to solve in seconds.  Times highlighted in red are faster than the corresponding method with small overlap.}
\begin{tabular}{|c|c|c|c|c|}
\hline
$\ell$ & $J=4$ & $J=16$ & $J=64$ & $J=256$ \\ \hline
 1 & \color{red}{1.3327$\times10^{-1}$} & - & - & -  \\ \hline
 2 & \color{red}{1.9001$\times10^{-1}$} & 1.6626$\times10^{+1}$ & - & -  \\ \hline
 3 & \color{red}{2.8851$\times10^{-1}$} & \color{red}{7.1764$\times10^{+0}$} & 6.1481$\times10^{+1}$ & -  \\ \hline
 4 & 5.8644$\times10^{+0}$ & \color{red}{2.5455$\times10^{+1}$} & \color{red}{4.5772$\times10^{+1}$} & 3.5521$\times10^{+2}$  \\ \hline
 5 & 1.0258$\times10^{+2}$ & \color{red}{7.3422$\times10^{+1}$} & \color{red}{6.8872$\times10^{+1}$} & \color{red}{1.5702$\times10^{+2}$}  \\ \hline
 6 & 1.3211$\times10^{+3}$ & 5.2160$\times10^{+2}$ & \color{red}{1.0986$\times10^{+2}$} & \color{red}{1.5081$\times10^{+2}$}  \\ \hline
 7 & 2.4185$\times10^{+4}$ & 8.1268$\times10^{+3}$ & \color{red}{7.7438$\times10^{+2}$} & \color{red}{3.0096$\times10^{+2}$}  \\ \hline
 8 & DNC & DNC & 1.1663$\times10^{+4}$ & \color{red}{1.6401$\times10^{+3}$}  \\ \hline
\end{tabular}
\label{tab:OLLO_t}
\end{table}
\begin{table}[htbp]
\caption{Two-level generous overlap: time to solve in seconds.  Times highlighted in red are  faster than the corresponding method with small overlap.}
\begin{tabular}{|c|c|c|c|c|}
\hline
$\ell$ & $J=4$ & $J=16$ & $J=64$ & $J=256$ \\ \hline
 1 & \color{red}{1.3327$\times10^{-1}$} & - & - & - \\ \hline
 2 & \color{red}{4.7170$\times10^{-1}$} & 1.6626$\times10^{+1}$ & - & - \\ \hline
 3 & \color{red}{6.4757$\times10^{-1}$} & \color{red}{1.0396$\times10^{+1}$} & 6.1481$\times10^{+1}$ & - \\ \hline
 4 & 6.7362$\times10^{+0}$ & \color{red}{3.4515$\times10^{+1}$} & \color{red}{6.3328$\times10^{+1}$} & 3.5521$\times10^{+2}$ \\ \hline
 5 & 1.0614$\times10^{+2}$ & \color{red}{8.1754$\times10^{+1}$} & \color{red}{8.7046$\times10^{+1}$} & 1.4895$\times10^{+2}$ \\ \hline
 6 & 1.3276$\times10^{+3}$ & 5.4679$\times10^{+2}$ & \color{red}{1.1548$\times10^{+2}$} & \color{red}{1.1275$\times10^{+2}$} \\ \hline
 7 & 2.3118$\times10^{+4}$ & 8.4160$\times10^{+3}$ & 7.5119$\times10^{+2}$ & \color{red}{1.9706$\times10^{+2}$} \\ \hline
 8 & DNC & DNC & 1.1962$\times10^{+4}$ & 1.1392$\times10^{+3}$ \\ \hline
\end{tabular}
\label{tab:TLLO_t}
\end{table}
\section{Conclusion}
 In this paper we present
 additive Schwarz preconditioners for
  the linear systems that arise from the PDAS algorithm
    applied to an elliptic distributed optimal control problem with pointwise state constraints
   discretized by a PUM.
 Based on the condition number estimates and the numerical results, the two-level method with small overlap appears to be the best choice
 for small $h$ and large $J$.
\section*{Acknowledgements}
The work of the first and third authors was supported in part
 by the National Science Foundation under Grant No.
 DMS-16-20273.
Portions of this research were conducted with high performance computing
resources provided by Louisiana State University (http://www.hpc.lsu.edu).

\end{document}